\documentclass[format=sigconf]{acmart}
\usepackage{fancyhdr}
\usepackage{graphicx}
\usepackage{amsmath}

\usepackage{amssymb}
\usepackage{mathtools}
\usepackage{comment}
\usepackage{subfig}
\usepackage{bm}
\usepackage{multirow}
\usepackage{threeparttable,booktabs}
\usepackage{tikz}
\usepackage{balance}
\usepackage{courier}                                       % courier font, used in \texttt
\usepackage[mathcal]{eucal}
\usepackage[]{algpseudocode}                               % algorithm package
\algrenewcommand\textproc{\texttt}
\usepackage{algorithm}
\usepackage{filecontents}                                  % support to pgfplots
\usepackage{pgfplots}
\usepackage{pgfplotstable}
\pgfplotsset{compat=newest}
\renewcommand\footnotetextcopyrightpermission[1]{} 

\theoremstyle{plain}
\newtheorem{Theorem}{\textbf{Theorem}}

\newtheorem{Proposition}{\textbf{Proposition}}
\theoremstyle{definition}

\usepackage[skip=1pt]{caption}            % set space between figure and caption
\makeatletter
\def\@mkbibcitation{\relax}
\makeatother
\begin{document}

\title{GBEM: Galerkin Boundary Element  Method for 3-D Capacitance Extraction}

% Single author syntax
\author{\mbox{Shengkun Wu$^{1}$} ~and~\mbox{Xingquan Li$^{2,1}$}}
\affiliation{\institution{$^{1}$ Peng Cheng National Laboratory, Shenzhen 518000, China. }}
\affiliation{\institution{$^{2}$ School of Mathematics and Statistics, Minnan Normal University, Zhangzhou 363000, China. }}
\email{shengkunwu@foxmail.com,  fzulxq@gmail.com}

%\date{}
\baselineskip=10.53pt

%%%%%%%%%%%%%%%%%%%%%%%%%%%%%%%%%%%%%%%%%%%%%%%%%%%%%%%%%%%%%%%%%%%%%%%%%%%%%%%%
% Abstract
%%%%%%%%%%%%%%%%%%%%%%%%%%%%%%%%%%%%%%%%%%%%%%%%%%%%%%%%%%%%%%%%%%%%%%%%%%%%%%%%
\begin{abstract}
For modern IC design, electromagnetic coupling among interconnect wires plays an increasingly important role in signoff analysis. The requirement of fast and accurate capacitance extraction is becoming more and more urgent.
The critical step of extracting capacitance among interconnect wires is solving electric field. However, due to the high computational complexity, solving electric field is extreme timing-consuming. To improve computational efficiency, we propose a Galerkin boundary element method (GBEM) to extract capacitance. The advantage of this method is that it can greatly reduce the number of boundary elements on the premise of ensuring that the error is small enough. As a consequence, the matrix order of the discretization equation will also decrease. The experiments in this paper have proved this advantage of our algorithm. Moreover, we have took advantage of some mathematical theorems in this paper. 
Our attempt shows that there will be more connections between the capacitance extraction and some mathematical conception so that we can use more mathematical tools to solve the problems of capacitance extraction.
\end{abstract}

 \vspace{-0.2cm}

\maketitle
%%%%%%%%%%%%%%%%%%%%%%%%%%%%%%%%%%%%%%%%
%% Introduction
%%%%%%%%%%%%%%%%%%%%%%%%%%%%%%%%%%%%%%%%
\section{Introduction}
\label{sec:introduction}

%介绍RCX背景，分2D和3D
With the feature size of integrated circuit scaling down, the coupling capacitance of interconnect wires is making more and more significant impact on circuit performance. The most of existed mature technologies of field solver are based on numerical approach.
There are mainly several capacitance extraction methods which are the finite difference method (FDM), the finite element method (FEM), the boundary element method (BEM) and the floating random walk (FRW) method. These methods are called field solver and their characteristics are summerized in Table \ref{table:compare}.

\begin{table}[h]
  \centering
  \footnotesize
  \caption{Features of different Field Solver.}
  \vspace{0.0cm}
  \begin{tabular}{|l|c|c|c|c|c|c|
  }
  \hline
      -- & FDM or FEM   & BEM    & FRW \\  
  \hline
   Equation form      & differential    & integral    & integral \\
   \hline
  Discretization mode   &  domain discretization  & boundary discretization   & --  \\
   \hline

  Generating Matrix   & large and sparse   & small and dense   & --  \\
  \hline
  Parallelism     & bad   & bad   & good  \\
   \hline
Convergence rate     & rapid   & rapid   & slow  \\
  
   \hline
   Main error source    & discretization   & discretization   & random  \\
   \hline
adaptability to &    &    &   \\
complex structure& good & bad & bad\\
   \hline
  \end{tabular}
  \label{table:compare}
\end{table}

The FEM and the FDM are classified as the domain discretization method.
It usually produces a sparse matrix with large order, see \cite{Zemanian1988}.  In 3-D capacitance extraction, because the order of the matrix increase rapidly, the speed of this method is limited. However, the domain discretization method is well established, thus this method is also used by some software. 

The FRW algorithm for capacitance extraction, presented as a 2-D version, was proposed in 1992 \cite{Coz} . Its idea is to convert the calculation of conductor charge to the Monte Carlo integration performed with FRWs.
The random walk method is advantageous in  parallelism over the traditional methods, see \cite{Yang2020}. Recently, the FRW method has been made a lot of progress. However, unlike the classical analytic method, the error of FRW is random.

%专门阐述BEM,优点：只需看边界；缺点：矩阵维度仍然较多
The boundary element method only needs to discretize the boundary, thus the matrix order produced by the BEM is smaller than that produced by the FDM. However, the matrix obtained by the BEM is not sparse and a lot of time is spent in calculating the matrix elements. An progress we will mention is the quasi-multiple medium (QMM) method \cite{Yu2003}.  By adding dielectric interfaces, the order of the coefficient matrix produced by the QMM method increases slightly but the matrix becomes sparse. In contrast, we will introduce our ideas to reduce the order of the matrix.

%引出我们的动机：希望找到一个降低矩阵维度的方法；希望从边界元的划分入手；旨在能够得到一个理论误差估计
The Garlerkin method for the Laplace equation has been studied in Mathematics for a long time, see \cite{Steinbach}. Its advantage is that it has an error estimation which can guide us to partition the boundary. With this guidance, we only need a small number of boundary elements to obtain high accuracy. However, the matrix elements obtained from the Garlerkin method is complicated. In fact,  we'll find that we have to deal with quadruple integrals. To make the algorithm faster, we will provide our suggestion about the calculation of the matrix elements.

The main contributions of the GBEM are summarized as follows:
\begin{itemize}
    %引入Galerkin 投影，得到对称矩阵A，Ax=b好求解
	\item We use the Galerkin method to improve the boundary element method of the capacitance extraction.
	
    %得到Q的理论误差评估界
	\item We achieve a theoretical error estimation for the single dielectric case. The error estimation is obtained by several theorems of the Galerkin method and integral operators involved in the boundary integral equation.
    
    %在理论误差的指导下，提出高效的边界元划分策略
	\item With the guide of our error estimation, we develop a boundary partition scheme. This boundary partition strategy can largely reduce the number of boundary elements and ensure sufficient accuracy.
    
    %实验结果说明我们的方法的有效性
	\item Numerical experiments show that, our approach can obtain accurate results with much fewer boundary elements.

\end{itemize}

The remainder of this paper is organized as follows. In Section \ref{sec:Preliminaries}, we provide the basic conception of the boundary element method. In Section \ref{single}, we propose the Galerkin boundary element method and provide our suggestion about calculating matrix elements for the single dielectric case.  In Section \ref{multi-dielectric}, we discuss the multi-dielectric case.
Section \ref{sec:experimental results} shows experimental results.

\section{Preliminaries} \label{sec:Preliminaries}

In this section, we formulate the capacitance extracting problem.

\subsection{Problem Statement} \label{sec:problem_statement}

\begin{figure}[h] \centering
  \includegraphics[width=8.7cm]{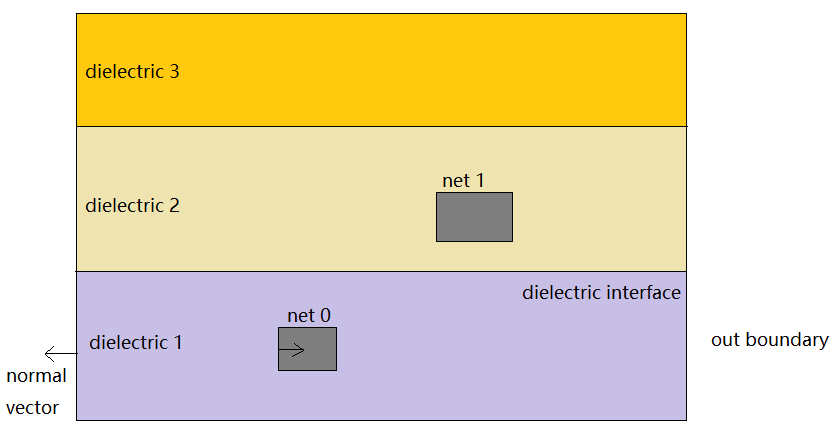}
  \caption{\small
      A cross section of 3-D capacitance extraction problem with multi-dielectrics
  }
  \label{fig:crosssection}
\end{figure}
In each dielectric region $\Omega_a$, the electric potential $u_a$ satisfies the Laplace equation
\begin{equation*}
\left\{
\begin{aligned}
& \nabla^2u_a=0 \text{ in } \Omega_a \\
& u_a=1 \text{ on the main net}\\
& u_a=0 \text{ on other nets and the Dirichlet boundary}\\
& \frac{\partial u_a}{\partial n}=0 \text{ on the Neumann boundary.}\\
\end{aligned}
\right.
\end{equation*}

Let $q_a$ and $q_b$ be the normal electric field intensity in regions $\Omega_a$ and $\Omega_b$ respectively. On the dielectric interface $\Gamma_{ab}$ of two dielectric region $\Omega_a$ and $\Omega_b$, the compatibility equation holds:

\begin{equation*}
\begin{aligned}
& \varepsilon_a q_a=-\varepsilon_b q_b \\ 
& u_a=u_b.
\end{aligned}
\end{equation*}
Using Green's identity, we have following integral equation:
\begin{equation}\label{representation}
	\sigma(x)u_a(x)+\int_{\partial \Omega_a} q^*(x,y) u_a(y)d\Gamma(y)
	=\int_{\partial \Omega_a} u^*(x,y)q_a(y)d\Gamma(y),
\end{equation}
where 
$$q^{*}(x,y)=\frac{\langle x-y,n_y \rangle}{4\pi |x-y|^3} \text{\quad, \quad  } 
u^*(x,y)=\frac{1}{4\pi |x-y|}$$
and $\sigma(x)$ satisfies 
$\sigma(x)=\frac{1}{2} \text{ for almost all } x\in \Gamma.$

\subsection{Classical Boundary Method}
The classical boundary element method employs boundary element partition and evaluating the boundary integral equation  at collocation points, one for an element. The discretized boundary integral equation in the region $\Omega_a$ is
\begin{equation*}
	\sigma(x) u(x)+\sum_{j}\int_{\Gamma_{j}} q^*(x,y)d\Gamma(y)u_j
	=\sum_{j}\int_{\Gamma_{j}} q(y)d\Gamma(y)q_j,
\end{equation*}
see \cite{Yu2014}. Combining this with the compatibility equation, the electric field intensity on each element can be obtained.

\section{single Dielectric Case} \label{single}
\subsection{Linear Equation}
In this section, we discuss the single dielectric case with the Dirichlet boundary condition
$$u=0 \quad \text{  on the outer boundary of the dielectric.}$$

Duo to the boundary conditions, we have
\begin{align*}
	\sigma(x)u(x)+\int_{-\partial mnet} q^*(x,y) d\Gamma(y)
	=\int_{\partial \Omega} u^*(x,y)q(y)d\Gamma(y),
\end{align*}
where $-\partial mnet$ represent the boundary of the main net and the normal vectors point to the internal side of the main net.

If $x$ is on the main net, we have 
\begin{align*}
	&\sigma(x)u(x)+\int_{-\partial mnet} q^*(x,y) d\Gamma(y)\\
	=&\frac{1}{2}-\int_{\partial mnet} q^*(x,y) d\Gamma(y)(\text{ for almost all } x \in \partial mnet).\\
\end{align*}
In the inner side of the main net, we suppose there is an electric potential $v$ which equals to 1 every where. Using (\ref{representation}), we obtain
\begin{equation*}
	\sigma(x)v(x)+\int_{\partial \Omega} q^*(x,y) v(y)d\Gamma(y)=0.
\end{equation*}
Thus
\begin{equation*}
	\frac{1}{2}+\int_{\partial \Omega} q^*(x,y)d\Gamma(y)=0 \quad (\text{ for almost all } x \in \partial mnet).
\end{equation*}
If $x$ is on other nets or the dielectric boundary, we have $u(x)=0$. Using Gauss's theorem, we get
$$\int_{-\partial mnet} q^*(x,y) d\Gamma(y)=0.$$

To summarize, we have
\begin{align*}
\int_{\partial \Omega} u^*(x,y)q(y)d\gamma(y)=f_0(x),
\end{align*}
where
\begin{equation*}
f_0(x)=\left\{
\begin{aligned}
&1 \text{\quad if x is on $\partial mnet$} \\
&0 \text{\quad otherwise.}
\end{aligned}
\right.
\end{equation*}

To solve $q$, we need to partition the boundary into finite pieces. 
We suppose $\partial \Omega=\bigcup_{j} I_j$ and use $q_j$ to approximate $q(x)$ when $x\in I_j$, then we get following equation:
\begin{align*}
\sum_{j}\int_{I_{j}} u^*(x,y)d\Gamma(y)q_{j}=f_0(x).
\end{align*}
Integrating both sides over $I_i$, we get
\begin{equation*}
\sum_{j}\int_{I_i}\int_{I_{j}} u^{*}(x,y)d\Gamma(y)\Gamma(x) q_{j} =\int_{I_i}f_0(x)d\Gamma(x)
\end{equation*}
and
\begin{equation}\label{myequation}
\sum_{j}\frac{1}{|I_i|}\frac{1}{|I_j|}\int_{I_i}\int_{I_{j}} u^{*}(x,y)d\Gamma(y)\Gamma(x) q_{j}|I_j| =\frac{1}{|I_i|}\int_{I_i}f_{0}(x)d\Gamma(x),
\end{equation}
where $|I_i|$ is the area of $I_i$. Let 
\begin{equation}\label{matrix}
A=\Big\{ \frac{1}{|I_i|}\frac{1}{|I_j|}\int_{I_i}\int_{I_{j}} u^{*}(x,y)d\Gamma(y)\Gamma(x) \Big\}_{i,j},
\end{equation}
$$X=\Big\{q_j|I_j|\Big\}_{j}\quad \text{and}\quad b=\Big\{ \frac{1}{|I_i|}\int_{I_i}f_{0}(x)d\Gamma(x)\Big\}_{i}.$$
Then, we obtained a linear equation
\begin{equation}\label{lineareq1}
Ax=b.
\end{equation}
Once we solved this equation, the electric charge on the net n we obtained is given by 
\begin{equation}\label{electriccharge}
Q^{'}_{n}=\sum_{I_j\subset \partial {netn}}\varepsilon q_j|I_j|, 
\end{equation}
where $\varepsilon$ is the dielectric constant.
To solve the equation (\ref{lineareq1}), we need to know more about the information of the matrix $A$. It is easy to see that $A$ is symmetric. 
In fact, we can show that $A$ is positive definite. 
Moreover, we will give an error estimation of the electric charge on each net, through this estimation we will state our strategy of how to partition the boundary.

\subsection{Error Estimation}
To achieve our goals, we need two theorems. Let $V$ be an operator such that
\begin{align*}
V q(x)=\int_{\partial \Omega} u^*(x,y)q(y)d\Gamma(y).
\end{align*}

\begin{Theorem}\cite[Theorem 6.22]{Steinbach}
\label{theorem:elliptic}
If the dimension of the region is three, then the operator $V$ is elliptic on the Sobolev space $H^{-1/2}(\partial\Omega)$  , i.e. 
\begin{align*}
	\langle V g,g \rangle \geq c \|g\|_{H^{-1/2}(\partial\Omega)}^2 \text{ for all $g$ in } H^{-1/2}(\partial\Omega).
\end{align*}
\end{Theorem}
By Theorem \ref{theorem:elliptic}, we know that the operator $V$ is positive definite, thus the matrix $A$ defined via (\ref{matrix}) is positive definite.
Once the matrix is symmetric and  positive definite, the Cholesky decomposition is an  efficient way to solve equation (\ref{lineareq1}).

\begin{Theorem}\cite[Lemma 8.1]{Steinbach}\label{theorem:Galerkin}
If $V$is a bounded operator from Banach space $X$ to its dual space $X^{'}$ and elliptic operator on $X$.
Let $q$ be a solution of the equation
\begin{equation*}V q=f \text{ given } f\in X^{'}.\end{equation*}
Let $X_{m}$ be a finite dimensional subspace of X, and $q^{'}\in X_{m}$ be the solution of following equation 
\begin{equation}\label{GalerkinEquation}
 \langle V q', p\rangle = \langle f,p \rangle \text{ for all } p\in X_{m}.
\end{equation}
We have
\begin{align*}
   \|q-q'\|_{X}\leq C \inf_{p\in X_m}\|q-p\|_{X}.
\end{align*}
\end{Theorem}

Equation (\ref{GalerkinEquation}) is called Galerkin Method. To apply Theorem \ref{theorem:Galerkin}, Let $X_m$ the linear span of functions $\Big\{\frac{\chi_{I_j}}{|I_{j}|}\Big\}$, where 
\begin{equation*}
\chi_{I_j}(x)=\left\{
\begin{aligned}
&1 \text{\quad if x is on $I_{j}$} \\
&0 \text{\quad otherwise.}
\end{aligned}
\right.
\end{equation*}
Thus $X_m$ is a finite dimensional subspace of $H^{-1/2}(\partial\Omega)$. Then (\ref{theorem:Galerkin}) is equivalent to
$$\langle V q', \frac{\chi_{I_i}}{|I_i|}\rangle 
= \langle f,\frac{\chi_{I_i}}{|I_i|} \rangle.$$
Since $q^{'}\in X_{m}$, we suppose $q^{'}=\sum_{j}c_j\chi_{I_j}$. Thus
\begin{align*}
\sum_{j}\frac{1}{|I_i|}\frac{1}{|I_j|}\int_{I_i}\int_{I_{j}} u^{*}(x,y)d\gamma(y)\Gamma(x) c_j|I_j| 
= \langle f,\frac{\chi_{I_i}}{|I_i|}\rangle.
\end{align*}

Comparing this equation with (\ref{myequation}), we know that equation (\ref{myequation}) coincides with the Galerkin method (\ref{GalerkinEquation}) and the relationship is given by
$$q^{'}=\sum_{j}q_j\chi_{I_j},$$
where $q^{'}$ is the solution of the Galerkin method (\ref{GalerkinEquation}) and $\{q_j\}$ comes from (\ref{myequation}).
We will use this relationship to provide an error estimation of our method.

\begin{Proposition}\label{error}
Let $Q_n$ be the electric charge on the net n. Let $Q^{'}_{n}$ be the electric charge on the net n given by (\ref{electriccharge}). The projection of $q\in H^{-1/2}(\partial(\Omega)$ into $X_{m}$ is given by $$P_m q=\sum_{j}\frac{1}{|I_j|}\int_{I_j}q(x)d\Gamma(x) \chi_{I_j},$$
then there is a constant such that 
$$|Q_n-Q^{'}_{n}|\leq C \|q-P_{m}q\|_{L^2(\partial \Omega)}.$$
If the tangential derivative of $q$ exists, we further have
$$|Q_n-Q^{'}_{n}|\leq C \sup_{j} 
\Big(|I_j|\sup_{x\in I_{j}}\frac{\partial q}{\partial \tau}(x)\Big) $$
\end{Proposition}
\textbf{Proof:}
Based on Gauss's Theorem, we have
\begin{align*}
|Q_n-Q^{'}_n|&=\Big|\varepsilon\int_{\partial netn}  q(x)d\Gamma(x)
-\sum_{I_j\subset \partial {netn}}\varepsilon q_j|I_j|\Big|\\
&=\varepsilon\Big|\int_{\partial netn}  q(x)-q^{'}(x)d\Gamma(x)
\Big|\\
&=\varepsilon\Big|\langle (q(x)-q^{'}(x)),\chi_{\partial netn}\rangle\Big|\\
&\leq \varepsilon\|q(x)-q^{'}(x)\|_{H^{-1/2}(\partial \Omega)}
*\|\chi_{\partial netn}\|_{H^{1/2}(\partial \Omega)}\\
&=\varepsilon|\partial netn| \cdot
\|q-q^{'}\|_{H^{-1/2}(\partial \Omega)}.
\end{align*}
Using Theorem \ref{theorem:Galerkin}, we obtain
\begin{align*}
 \varepsilon|\partial netn|\cdot
\|q-q^{'}\|_{H^{-1/2}(\partial\Omega)}
\leq C\varepsilon|\partial netn| \inf_{p\in X_{m}}\|q-p\|_{H^{-1/2}(\partial\Omega)}.
\end{align*}
Thus 
\begin{align*}
|Q_n-Q^{'}_n| &\leq C\varepsilon|\partial netn| \inf_{p\in X_{m}}\|q-p\|_{H^{-1/2}(\partial\Omega)}\\
&\leq C\varepsilon|\partial netn| \inf_{p\in X_{m}}\|q-p\|_{L^2(\partial \Omega)}\\
&= C\varepsilon|\partial netn|\|q-P_{m}q\|_{L^2(\partial\Omega)}.
\end{align*}
Further, if the tangential derivative of $q$ exists, then
\begin{align*}
&|Q_n-Q^{'}_n|\\
\leq& C\varepsilon|\partial netn| 
\Big(\int_{\partial\Omega} \Big|q(x)-\sum_{j}\frac{1}{|I_j|}\int_{I_j}q(y)d\Gamma(y)\chi_{I_j}(x)\Big|^{2}d\Gamma(x) \Big)^{1/2}\\
=& C\varepsilon|\partial netn| 
\Big(\int_{\partial\Omega} \Big|\sum_{j}\frac{1}{|I_j|}\int_{I_j}q(x)-q(y)d\Gamma(y)\chi_{I_j}(x)\Big|^{2}d\Gamma(x) \Big)^{1/2}\\
=& C\varepsilon|\partial netn| 
\Big(\int_{\partial\Omega} \Big|\sum_{j}\frac{1}{|I_j|}\int_{I_j}|I_j|\sup_{x\in I_{j}}|\frac{\partial q}{\partial \tau}(x)|d\Gamma(y)\chi_{I_j}(x)\Big|^{2}d\Gamma(x) \Big)^{1/2}\\
\leq & C\varepsilon|\partial netn| 
|\partial \Omega|^{1/2}
\sup_{j}\Big(|I_j|\sup_{x\in I_{j}}|\frac{\partial q}{\partial \tau}(x)|\Big)^{1/2},
\end{align*}
to complete the proof.
$\hfill{} \Box$

Proposition \ref{error} tells us the error of the electric charge on each net. Since the capacitance of the main net is given by $C_{1}=Q_{1}$ and the coupling capacitance between the main net and the net n is given by  $C_{n}=Q_{n}$, we obtained the error estimation of each capacitance.

To understand this proposition, we need some discussion. Proposition \ref{error} tells us that the error is dominated by  the area of the boundary element times the derivative of the electric field intensity on that element. On the other hand, the boundary conditions tell us that the main net has electric potential 1 and other boundaries have electric potential 0, thus the derivative of the electric field intensity decreases rapidly as the distance to the main net increases. Thus, for the boundary element which is far away from the main net, the area of the boundary element can be very large due to the derivative of the electric potential is small. As a result, for the boundary element which is far away from the main net, we only need to partition it into several pieces. That is the reason that this method can reduce the number of boundary elements. Next, we will provide several formulas for the boundary partition strategy.

\subsection{Partition Boundary}\label{partition}
\begin{figure}[h] \centering
  \includegraphics[width=8.7cm]{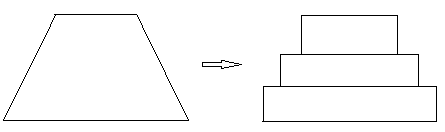}
  \caption{\small
      Using cuboid to approximate trapezoid.
  }
  \label{trapezoid}
\end{figure}

We suppose all conductors are construct by cuboid. If a conductor is a trapezoid, we use several cuboids to approximate it, see figure \ref{trapezoid}.

Our boundary element partition strategy is based on the Proposition \ref{error}.
If we want to make the error small, we need to partition $\partial\Omega$ such that
$$|I_j|\sup_{x\in I_{j}}|\frac{\partial q}{\partial \tau}(x)|\leq c$$
for some constant $c$ which is small. That is to say we want
$$|I_j|\leq\frac{c}{\sup_{x\in I_{j}}|\frac{\partial q}{\partial \tau}(x)|}.$$

Because the region is irregular, it's very hard to obtain priori estimates of $|\frac{\partial q}{\partial \tau}(x)|$. However, we can take a physical approximation from the boundary conditions of the capacitance extraction problem.

According to the boundary conditions, the electric potential of the main net is 1 and the electric potential of the rest boundary is 0. Thus, the electric field intensity decreases as the distance to the main net increases. We can imagine that if the boundary is far away from the main net, then $|\frac{\partial q}{\partial \tau}(x)|$ is small, as a consequence the area of the responding boundary element can be large. If $x$ is far away from the main net, then we can regard the main net as a ball with radio $r$ and center $x_0$. If the electric potential of the ball is 1, then the electric field density of point $x$ is given by 
$$q(x)=\frac{r^2}{|x-x_0|^2}\cdot \frac{x-x_0}{|x-x_0|}.$$
Thus we have
\begin{equation}\label{electric_derivative}
|\frac{\partial q}{\partial \tau}(x)|\leq \frac{C}{|x-x_0|^{3}}.
\end{equation}
Our partition strategy for the boundary which is far away from the main net is based on equation (\ref{electric_derivative}).

Now, we propose our strategy of partition. 
Let $p_1,p_2,p_3,p_4,p_5$ are parameters which can be choose for different cases .

If the distance between the boundary and the main net is less than $p_5$, then we partition it into rectangles $\{I_j\}$ such that 
$$|I_j|\leq p_4.$$

If the distance of the boundary is larger than $p_5$, then we partition it into rectangles $\{I_j\}$ such that 
\begin{equation*}
|I_j|\leq p_4*p_3*(d_{I_j}/p_5)^{3} ,
\end{equation*}
where $d_{I_j}$ is the distance between $I_j$ and the main net.
Duo to the electric field shielding effect, if the boundary is back to the main net then we set
$p_4=p_1\cdot p_2.$
If the boundary faces to the main net, we set 
$p_4=p_1$. Here $p_2$ is usually larger than 1.

\begin{figure}[htbp] \centering
  \includegraphics[width=8.7cm]{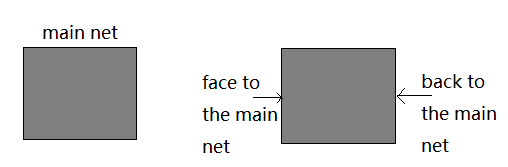}
  \caption{\small
      Different boundary location
  }
  \label{fig:back}
\end{figure}

\subsection{Matrix Element}\label{uintegral}
The Galerkin method has brought us both a benefit and a trouble. The benefit is that the number of boundary elements is reduced. The trouble is that the matrix elements are complicated.
The matrix elements are given by
\begin{align*} a_{i,j}=& \frac{1}{|I_i|}\frac{1}{|I_j|}\int_{I_i}\int_{I_{j}} u^{*}(x,y)d\Gamma(y)\Gamma(x)\\
=&\frac{1}{4\pi|I_i||I_j|}\int_{I_i}\int_{I_{j}} \frac{1}{|x-y|}d\Gamma(y)\Gamma(x).
\end{align*}
Since $I_{i}$ and $I_{j}$ are rectangles, this integral is actually a quadruple integral. Directly using numerical integration to calculate $a_{i,j}$ will take a lot of time. Thus we have to use some tricks. 

To handle the integral, we need a notation. If $f(x)$ is a function, we denote
$$D_{x}(x_{u}, x_{d})f(x)=f(x_{u})-f(x_{d})$$
and
\begin{align*}
&D_{x_{1},\cdots,x_n}(x_{1u}, x_{1d},\cdots,x_{nu}, x_{nd})\\
=&D_{x_{1}}(x_{1u}, x_{1d})D_{x_{2},\cdots,x_n}(x_{2u}, x_{2d},\cdots,x_{nu}, x_{nd}).
\end{align*}
Let $\vec{n}(I_i)$ and $\vec{n}(I_j)$ be the normal vector of $I_i$ and $I_j$ respectively.

If $\vec{n}(I_i)$  is parallel to $\vec{n}(I_j)$, suppose
\begin{equation}\label{parallel}
\begin{split}
&I_i=\{(x_1,x_2,x_3): x_{1d}\leq x_1\leq x_{1u}, x_{2d}\leq x_2\leq x_{2u}\}\\
&I_j=\{(y_1,y_2,y_3): y_{1d}\leq y_1\leq y_{1u}, y_{2d}\leq y_2\leq y_{2u}\}.
\end{split}
\end{equation}
Let $b=x_3-y_3$. Then 
\begin{align*}
&\int_{I_i}\int_{I_{j}} \frac{1}{|x-y|}d\Gamma(y)\Gamma(x)\\
=&\int_{x_{1d}}^{x_{1u}}\int_{x_{2d}}^{x_{2u}}\int_{y_{1d}}^{y_{1u}}\int_{y_{2d}}^{y_{2u}} \frac{1}{\sqrt{(x_{1}-y_{1})^2+(x_{2}-y_{2})^2+b^2}}dy_{2}dy_{1}
dx_{2}dx_{1}\\
=&\int_{x_{1d}}^{x_{1u}}\int_{x_{2d}}^{x_{2u}}\int_{-y_{1d}}^{-y_{1u}}\int_{-y_{2d}}^{-y_{2u}} \frac{1}{\sqrt{(x_{1}+y_{1})^2+(x_{2}+y_{2})^2+b^2}}dy_{2}dy_{1}
dx_{2}dx_{1}\\
=&D_{x_1,y_1}(x_{1u},x_{1d},-y_{1u},-y_{2d})\Bigg[\int_{x_{2d}}^{x_{2u}}
\int_{-y_{2d}}^{-y_{2u}}\\
&(x_{1}+y_{1})\ln\big[(x_{1}+y_{1})+\sqrt{(x_{1}+y_{1})^2+(x_{2}+y_{2})^2+b^2}\big]dy_{2}dx_{2}\\
&-\int_{x_{2d}}^{x_{2u}}\int_{-y_{2d}}^{-y_{2u}} \sqrt{(x_{1}+y_{1})^2+(x_{2}+y_{2})^2+b^2}dy_{2}dx_{2}\Bigg].
\end{align*}
On one hand we have
\begin{equation}\label{int1}
\begin{split}
&\int_{x_{2d}}^{x_{2u}}
\int_{-y_{2d}}^{-y_{2u}}
a\ln\big[a+\sqrt{a^2+(x_{2}+y_{2})^2}\big]dy_{2}dx_{2}\\
=& D_{x_{2},y_{2}}(x_{2u},x_{2d},-y_{2u},-y_{2d})G_{1}(x_{2}+y_{2}),
\end{split}
\end{equation}
where
\begin{align*}
G_{1}(u)=&\frac{a}{2}u^2\ln[a+\sqrt{a^2+u^2}]-\frac{3au^2}{4}-\frac{a^2}{2}\sqrt{u^2+a^2}\\
&+a^2u\ln(u+\sqrt{u^2+a^2}).
\end{align*}
On the other hand
\begin{equation}\label{int2}
\begin{split}
&\int_{x_{2d}}^{x_{2u}}\int_{-y_{2d}}^{-y_{2u}} \sqrt{a^2+(x_{2}+y_{2})^2}dy_{2}dx_{2}\\
=& D_{x_{2},y_{2}}(x_{2u},x_{2d},-y_{2u},-y_{2d})G_{2}(x_{2}+y_{2}),
\end{split}
\end{equation}
where
\begin{align*}
G_{2}(u)=&\frac{1}{6}[u^2+a^2]^{3/2}+\frac{a^2}{2}u\ln(u+\sqrt{u^2+a^2})-\frac{a^2}{2}\sqrt{u^2+a^2}.
\end{align*}
Using (\ref{int1}) and (\ref{int2}) we can obtain the analytic expression of $a_{i,j}$ in the case of $\vec{n}(I_i)$  is parallel to $\vec{n}(I_j)$ and $x3-y3=0$.

If $x3-y3\neq 0$, we can't obtain an analytic expression of $a_{i,j}$. Thus we have to use numerical integral. To speed up, we need to transform the double integral into a single integral.
If $G(u)$ is a function, then
\begin{equation}\label{int3}
\begin{split}
&\int_{x_{d}}^{x_{u}}\int_{y_{d}}^{y_{u}}G(x+y)dydx\\
=&\int_{0}^{\delta x}\int_{0}^{\delta y}G(x+y+x_{d}+y_{d})dydx \Big( \delta x=x_{u}-x_{d},\delta y=y_{u}-y_{d} \Big)\\
=&\int_{0}^{\delta y} uG(u+x_{d}+y_{d})du+\int_{\delta x}^{\delta y} \delta yG(u+x_{d}+y_{d}) du\\
 &+\int_{\delta x}^{\delta y+\delta x} (\delta y+\delta x-u)G(u+x_{d}+y_{d})du.
\end{split}
\end{equation}
Using (\ref{int3}), we can calculate $a_{i,j}$ by numerical integral.

If $\vec{n}(I_i)$  is orthogonal to $\vec{n}(I_j)$, we suppose
\begin{equation}\label{orthogonal}
\begin{split}
&I_i=\{(x_1,x_2,x_3): x_{1d}\leq x_1\leq x_{1u}, x_{3d}\leq x_3\leq x_{3u}\}\\
&I_j=\{(y_1,y_2,y_3): y_{1d}\leq y_1\leq y_{1u}, y_{2d}\leq y_2\leq y_{2u}\}.
\end{split}
\end{equation}
Let $\tilde{y}_{2u}=y_{2u}-x_{2}$, $\tilde{y}_{2d}=y_{2d}-x_{2}$,
$\tilde{x}_{3u}=x_{3u}-y_{3}$ and $\tilde{x}_{3d}=x_{3d}-y_{3}$, then
\begin{align*}
&\int_{I_i}\int_{I_{j}} \frac{1}{|x-y|}d\Gamma(y)\Gamma(x)\\
=&\int_{x_{1d}}^{x_{1u}}\int_{\tilde{x}_{3d}}^{\tilde{x}_{3u}}\int_{y_{1d}}^{y_{1u}}\int_{\tilde{y}_{2d}}^{\tilde{y}_{2u}} \frac{1}{\sqrt{(x_{1}-y_{1})^2+y_{2}^2+x_{3}^2}}dy_{2}dy_{1}
dx_{3}dx_{1}\\
=&D_{x1,y1}(x_{1u},x_{1d},-y_{1u},-y_{1d})\Big(\\
&\int_{\tilde{x}_{3d}}^{\tilde{x}_{3u}}\int_{\tilde{y}_{2d}}^{\tilde{y}_{2u}} \sqrt{y_{2}^2 +x_3^2+(x_1+y_1)^2}dy_2dx_3-\\
&\int_{\tilde{x}_{3d}}^{\tilde{x}_{3u}}\int_{\tilde{y}_{2d}}^{\tilde{y}_{2u}}(x_1+y_1)\ln\big[(x_1+y_1)+\sqrt{(x_1+y_1)^2+y_2^2+x_3^2}\big]dy_2dx_3 \Big).
\end{align*}
On one hand
\begin{equation}\label{int4}
\begin{split}
\int_{x_d}^{x_{u}}\int_{y_{d}}^{y_{u}} \sqrt{y^2 +x^2+a^2}dydx
=\int_{x_d}^{x_{u}}D_{y}(y_u,y_d)G_3(x,y) dx,
\end{split}
\end{equation}
where 
\begin{align*}
G_{3}(x,y)=\frac{y}{2}\sqrt{y^2+x^2+a^2}+\frac{x^2+a^2}{2}\ln\big[y+\sqrt{x^2+y^2+a^2}\big].
\end{align*}
On the other hand, we have
\begin{equation}\label{int5}
\begin{split}
&\int_{x_d}^{x_{u}}\int_{y_{d}}^{y_{u}}a\ln\big[a+\sqrt{a^2+y^2+x^2}\big]dydx\\
=&D_{y}(y_u,y_d)\int_{x_d}^{x_{u}}\frac{ya}{2}\ln\big[a+\sqrt{a^2+y^2+x^2}\big]dx\\
&+D_{x}(x_u,x_d)\int_{y_d}^{y_{u}}\frac{xa}{2}\ln\big[a+\sqrt{a^2+y^2+x^2}\big]dy\\
&+\frac{1}{2}\int_{x_d}^{x_{u}}\int_{y_{d}}^{y_{u}}\Big(a-\frac{a^2}{\sqrt{a^2+y^2+x^2}}\Big)dydx,
\end{split}
\end{equation}
where the last double integral is given by
\begin{align*}
&\int_{x_d}^{x_{u}}\int_{y_{d}}^{y_{u}}\Big(a-\frac{a^2}{\sqrt{a^2+y^2+x^2}}\Big)dydx\\
=&D_{x,y}(x_u,x_d,y_u,y_d)\Big[axy\Big]\\
&-D_{y}(y_u,y_d)\int_{x_d}^{x_{u}}a^2\ln(y+\sqrt{y^2+x^2+a^2}) dx.
\end{align*}

To summarize, we can calculate $a_{i,j}$ by several single numerical integrals or analytic expression. We used the Romberg quadrature formula to calculate each numerical integral in our program.

\section{multi-dielectric case} \label{multi-dielectric}
\subsection{Linear Equation}
In this section, we discuss the multi-dielectric case. Duo to equation (\ref{representation}), we have
\begin{equation*}
	\sigma(x)u_a(x)+\int_{\partial \Omega_a} q^*(x,y) u_a(y)d\Gamma(y)
	=\int_{\partial \Omega_a} u^*(x,y)q_a(y)d\Gamma(y),
\end{equation*}
holds for $x\in\Omega_a$.
We partition the boundary such that 
$\partial \Omega=\bigcup_{j}I_{j}.$
Taking integral over $I_i$ on both side and approximating $q_a (y)$ and $u_a(y)$ on each element with constant, we then obtain
\begin{align*}
&\frac{1}{2} \int_{I_i}u_{a}d\Gamma(x)
+\sum_{k}\int_{I_i}\int_{I_{k}} q^*(x,y) d\Gamma(y)d\Gamma(x)u_{a,k}\\
=&\sum_{j}\int_{I_i}\int_{I_{j}} u^*(x,y)d\Gamma(y)d\Gamma(x) q_{a,j}
\end{align*}

We have obtained a linear equation. In the multi-dielectric case the operators in the equation above may not be elliptic operators, thus we don't have an error estimation. In fact, this problem can be solved via elliptic operators and that method can also give us an error estimation. However, the computational complexity of that algorithm is too high, so we don't use that algorithm.

We still use the partition strategy in Section \ref{partition} for the boundary element on conductors and out boundary. But, we need to  change the partition strategy for the boundary elements in the dielectric interface. Since the unknowns corresponding to these boundary elements include electric field intensity and electric field. If $x$ is far away from the main net, we suppose the electric potential satisfies
$$\Big|\frac{\partial u}{\partial \tau}(x)\Big|\leq \frac{C}{|x-x_0|^2}.$$
Thus, for the dielectric interface, our partition strategy is as follows.
If the distance between the boundary and the main net is less that $p_5$, then we partition it into rectangles $\{I_j\}$ such that 
$$|I_j|\leq p_4.$$
If the distance of the boundary is larger that $p_5$, then we partition it into rectangles $\{I_j\}$ such that 
$$|I_j|\leq p_4*p_3*(d_{I_j}/l)^{2},$$
where $d_{I_j}$ is the distance between $I_j$ and the main net.
\subsection{Matrix Elements}
The matrix elements include two kind of integrals:
$$\int_{I_i}\int_{I_{j}} u^*(x,y)d\Gamma(y)d\Gamma(x)\text{\quad and\quad} Q_{i,j}=\int_{I_i}\int_{I_{j}} q^*(x,y)d\Gamma(y)d\Gamma(x).$$
The first one have been discussed in Section \ref{uintegral}. We need to discuss the second one.

If $\vec{n}(I_i)$  is parallel to $\vec{n}(I_j)$, we suppose $I_i$ and $I_j$ satisfies (\ref{parallel}). We will use same symbols in Section \ref{uintegral}. 

If $b=x3-y3=0$, then by the formular of $q^*(x,y)$ we have $Q_{i,j}=0$.

If $b=x3-y3\neq 0$, similar to Section \ref{uintegral}, we have
\begin{align*}
Q_{i,j}=&\int_{x_{1d}}^{x_{1u}}\int_{x_{2d}}^{x_{2u}}\int_{y_{1d}}^{y_{1u}}\int_{y_{2d}}^{y_{2u}}\frac{\langle x-y,n_y \rangle}{4\pi |x-y|^3}dy_{2}dy_{1}dx_{2}dx_{1}\\
=&D_{x_1,y_1}(x_{1u},x_{1d},-y_{1u},-y_{2d})\int_{x_{2d}}^{x_{2u}}\int_{-y_{2d}}^{-y_{2u}}\\
&\frac{1}{4\pi}\frac{b\sqrt{(x_2+y_2)^2+(x_1+y_1)^2+b^2}}{(x_2+y_2)^2+b^2}dx_2dy_2.
\end{align*}
using (\ref{int3}), we can calculate $Q_{i,j}$.

If $\vec{n}(I_i)$  is orthogonal to $\vec{n}(I_j)$, we suppose $I_i$ and $I_j$ satisfies (\ref{orthogonal}). Then
\begin{align*}
    &Q_{i,j}\\
    =&\int_{x_{1d}}^{x_{1u}}\int_{x_{3d}}^{x_{3u}}\int_{y_{1d}}^{y_{1u}}\int_{y_{2d}}^{y_{2u}}\frac{x_3-y_3}{4\pi |x-y|^3}dy_{2}dy_{1}dx_{3}dx_{1}\\
    =&D_{x_3}(\tilde{x}_{3u},\tilde{x}_{3d})\int_{x_{1d}}^{x_{1u}}\int_{\tilde{y}_{1d}}^{y_{1u}}\int_{\tilde{y}_{2d}}^{y_{2u}}\frac{\sqrt{(x_1-y_1)^2+y_2^2+x_3^2}}{4\pi}dy_{2}dy_{1}dx_{1}\\
    =&D_{x_3,x_1,y_1}(\tilde{x}_{3u},\tilde{x}_{3d},x_{1u},x_{1d},-y_{1u},-y_{1d})\int_{\tilde{y}_{2d}}^{\tilde{y}_{2u}}\\
    &(x_1+y_1)\ln\Big[(x_1+y_1)+\sqrt{(x_1+y_1)^2+y_2^2+x_3^2}\Big]\\
    &-\sqrt{(x_1+y_1)^2+y_2^2+x_3^2}dy_2.
\end{align*}
Thus, we can calculate $Q_{i,j}$ by using several single numerical integrals.

\section{Experimental Results} \label{sec:experimental results}
In this section, we show two experiments.
Our experiments are carried on a Inter(R) Xeon(R) Gold 5218 server with CPU at 2.3 GHz.

\begin{figure}[h] \centering
  \includegraphics[width=5 cm]{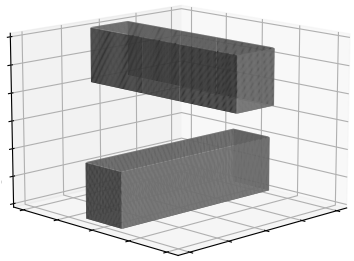}
  \caption{\small
      Single-dielectric test case.
  }
  \label{fig:single}
\end{figure}

We first present a single dielectric case. Two conductors with weight 1, height 1 and length 4 are immersed in a 40*40*10 region.  The lower left coordinate and the upper right coordinate of the net 1 are (-2, -0.5, 1) and (2, 0.5, 2). The lower left coordinate and the upper right coordinate of the net 2 are (-0.5, -2, 3) and (0.5, 2, 4). The lower left coordinate and the upper right coordinate of the dielectric boundary are (-20, -20, 0) and (20, 20, 10). The unit for each number is $\mu m$. We suppose the electric potential is 0 on the dielectric boundary.

\begin{table}[h]
  \centering
  \footnotesize
  \caption{Capacitances calculated with GBEM (in unit of $10^{-18}$ F).}
  \vspace{0.0cm}
  \begin{tabular}{|l|c|c|c|c|c|c|
  }
  \hline
      -- & $C_{11}$  &  $C_{12}$   & $C_{21}$    & $C_{22}$ & Time(s)& Memory(MB) \\  
  \hline
   74 panels   & 222   & -58.6   & -58.6   & 178 & 0.06&  7.19\\
   \hline
   183 panels   & 225   & -58.3   & -58.3   & 179 & 0.16 &  7.31 \\
  \hline
   594 panels   & 226   & -58.7   & -58.7   & 179  & 1 &12 \\
   \hline
  \end{tabular}
  \label{table:mysingle}
\end{table}

The result is presented in table \ref{table:mysingle}.
The result obtain by \cite{Sun1997} is $C_{11}=230$, $C_{12}=-61$, $C_{21}=-61$ and $C_{22}=180.6$, where the unit is $10^{-18}$ F. Thus, 74 boundary elements can give results with an error of less than 5\%.

Our multi-dielectric test case with six conductors comes from \cite[Page 56]{ Yu2014}. In this case the top and bottom dielectric boundary are the Dirichlet boundary. Other dielectric boundary are the Neumann boundary. Other data of this case can be found in \cite[Page 56]{ Yu2014}.

\begin{figure}[h] \centering
  \includegraphics[width=8.7cm]{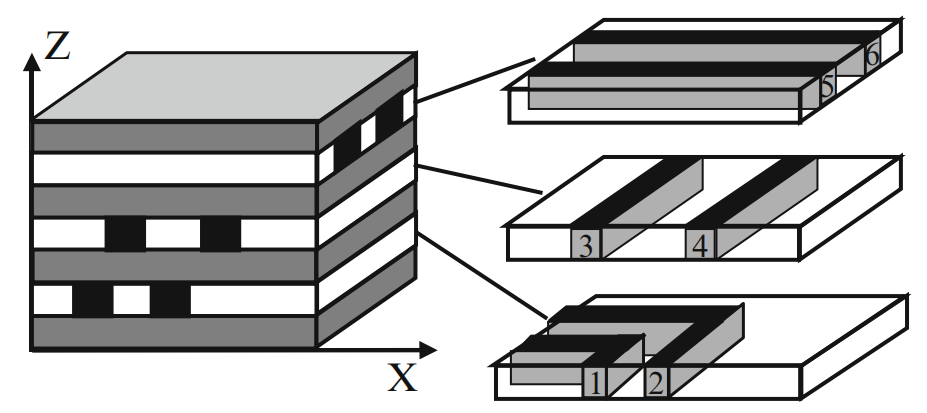}
  \caption{\small
      Multi-dielectric test case.
  }
  \label{fig:multi}
\end{figure}

\begin{table}[h]
  \centering
  \footnotesize
  \caption{The digonal entries of the capacitance matrices calculated with QMM (in unit of $10^{-15}$ F).}
  \vspace{0.0cm}
  \begin{tabular}{|l|c|c|c|c|c|c|
  }
  \hline
      -- & $C_{11}$  &  $C_{22}$   & $C_{33}$    & $C_{44}$   & $C_{55}$ & $C_{66}$ \\  \hline
   Capacitance &  0.682  &  1.31   & 1.6    & 1.54   & 2.53 & 2.53 \\ 
  \hline
  \end{tabular}
  \label{table:1}
\end{table}

In table \ref{table:1}, we provide the capacitance obtain by \cite{Yu2014} using QMM.
Because we used different server, we don't compare the run time. However, there are boundary elements from 2277 to 2575 in six computations for capacitance with QMM. We will show that, our method only need boundary elements from 135 to 316 to obtain the results such that the discrepancy of each number is less than 5\%.

\begin{table}[h]
  \centering
  \footnotesize
  \caption{The digonal entries of the capacitance matrices calculated with GBEM (in unit of $10^{-15}$ F).}
  \vspace{0.0cm}
  \begin{tabular}{|l|c|c|c|c|c|c|
  }
  \hline
      -- & $C_{11}$  &  $C_{22}$   & $C_{33}$    & $C_{44}$   & $C_{55}$ & $C_{66}$ \\  \hline
   Panel   & 135   & 159   & 159   & 164   & 165   & 170 \\
  \hline
   Time(s) & 0.165 & 0.203 & 0.186 & 0.156 & 0.156 & 0.158 \\ 
   \hline
   Memory(MB) &   7.36  &  7.36   & 7.39    & 7.49   & 7.35 & 7.33 \\ 
  \hline
   Capacitance &  0.685  &  1.307   & 1.65    & 1.58   & 2.56 & 2.56 \\
  \hline\\
  \hline
  Panel &  164  &  215   & 241    & 242   & 208 & 213 \\
  \hline
   Time(s) &   0.205  &   0.3   & 0.24    & 0.242   & 0.213 & 0.224 \\ 
   \hline
   Memory(MB) &   7.34  &  7.58   & 7.84    & 7.82   & 7.15 & 7.14 \\ 
  \hline
   Capacitance &  0.694  &  1.317   & 1.613    & 1.545   & 2.56 & 2.56 \\ 
  \hline\\
  \hline
  Panel &  217  &  316   & 312    & 314   & 297 & 306 \\
  \hline
   Time(s) &   0.251  &   0.422   & 0.394    & 0.402   & 0.375 & 0.371 \\ 
   \hline
   Memory(MB) &   7.26  &   8.375   & 8.38    & 8.47   & 8.17 & 8.16 \\ 
  \hline
   Capacitance &  0.691  &  1.306   & 1.598    & 1.542   & 2.543 & 2.539 \\ 
  \hline
  \end{tabular}
  \label{table:3}
\end{table}
In Table \ref{table:3}, we present our result.  For example, for $C_{11}$, if  the number of the boundary elements is 135, the run time would be 0.165s and the capacitance we obtained is $C_{11}=0.685$. We can find that the numbers of boundary elements required for our method are largely less than that for QMM. 

The details of our partition strategy is as follows. To obtain $C_{11}=0.691$, there are 217 boundary elements. In this situation, the main net is the net 1.
For each boundary $I$, let $d_{I}(\mu m)$ be the distance between $I$ with the net 1.
If $d_{I}$ is less than 1($\mu$m), we partition it in to $\cup_{j} I_{j}$ such that 
$$|I_{j}| < 1.5 \mu m^2. $$
For the Dirichlet boundary and the Neumann boundary, we partition it to $\cup_{j} I_{j}$ such that
$$|I_{j}| < p* 1.5 * d_{I_j}^3(\mu m^2) \text{ and } |I_{j}| < 1.5 * d_{I_j}^2(\mu m^2) \text{ respectively},$$
where $p=4$ if the boundary back to the main net and $p=1$ for other cases.

\section{Concluding Remarks} \label{sec:conclusion}

We proposed the Garlerkin boundary element method for the capacitance extraction. Based on the error analysis of this method, we proposed a boundary element partition strategy which can largely reduce the number of boundary elements. To calculate the matrix element quickly, we proposed a  method to deal with these integrals.

Because we don't have an algorithm which have both theoretical support and good experimental effect for the multi-dielectric case, this case is still need to be studied. What we have done is to provide an algorithm which have a good experimental results after many tests. In fact, the problems in multi-dielectric situations is call transmission problems in mathematics, see \cite{Gwinner}. The transmission problems is about the transmission of different fields in different regions. Capacitance extraction problems is a special case which only cares about the transmission of electric field in chip. on the one hand, it is possible to find out more mathematical tools to improve our algorithm. On the other hand, the capacitance extraction problems in IC design can also provide new challege to mathematics.

\end{document}